\documentclass[secthm,aap,dvips]{arximspdf}
\usepackage{mathbh}
\usepackage{graphicx}

\doi{10.1214/105051607000000168}
\volume{17}
\issue{4}
\pubyear{2007}
\firstpage{1347}
\lastpage{1361}

\makeatletter
\newtheorem{lem}{Lemma}[section]
\newproclaim{defn}{Definition}[section]
\newproclaim{rem}{Remark}[section]
\newproclaim{example}{Example}[section]

\makeatother

\begin{document}
\begin{frontmatter}

\title{A central limit theorem for stochastic recursive sequences of
topical operators}
\runtitle{CLT for SRS of topical operators}

\begin{aug}
\author[A]{\fnms{Glenn} \snm{Merlet}\corref{}\ead[label=e1]{glenn.merlet@gmail.com}}
\runauthor{G. Merlet}
\affiliation{CEREMADE\textup{,} Universit\'e Paris--Dauphine}
\address[A]{CEREMADE\\
Universit\'e Paris-Dauphine\\
Place du Mar\'echal de Lattre de Tassigny\\75775 Paris cedex 16\\
France\\
\printead{e1}} %adresu isvedimo komanda gale!
\end{aug}

% HISTORY:
\received{\smonth{8} \syear{2006}}
\revised{\smonth{1} \syear{2007}}

% ABSTRACT
%
\begin{abstract}
Let $(A_n)_{n\in\mathbb{N}}$ be a stationary sequence of topical
(i.e., isotone and additively homogeneous) operators. Let $x(n,x_0)$
be defined by $x(0,x_0)=x_0$ and $x(n+1,x_0)=A_nx(n,x_0)$.
It can model a wide range of systems including train or queuing
networks, job-shop, timed digital circuits or parallel
processing systems.

When $(A_n)_{n\in\mathbb{N}}$ has the memory loss property,
$(x(n,x_0))_{n\in\mathbb{N}}$ satisfies a strong law of
large numbers. We show that it also satisfies the CLT if $(A_n)_{n\in
\mathbb{N}}$
fulfills the same mixing and integrability assumptions that
ensure the CLT for a sum of real variables in the
results by P.~Billingsley and I.~Ibragimov.
\end{abstract}

% KEYWORDS
%

\end{frontmatter}

%s1 ###
\section{Model}

An operator $A\dvtx\mathbb{R}^d\rightarrow\mathbb{R}^d$ is called
additively homogeneous if it satisfies $A(x+a\mathbf{1})=Ax+a\mathbf
{1}$ for all $x\in\mathbb{R}^d$ and $a\in\mathbb{R}$, where
$\mathbf{1}$ is the vector $(1,\ldots,1)'$ in $\mathbb{R}^d$. It is
called isotone if $x\le y$ implies $Ax\le Ay$, in which the order is
the product order on $\mathbb{R}^d$. It is called topical if it is
isotone and homogeneous. The set of topical operators on
$\mathbb{R}^d$ will be denoted by~$\mathit{Top}_d$.

We recall that the action of matrices with entries in the semiring
${\mathbb R}_{\max}=(\mathbb{R}\cup\{-\infty\},\max,+)$ on
${\mathbb R}_{\max}^d$ is defined by $(Ax)_i=\max_j(A_{ij}+x_j)$.
When matrix $A$ has no $-\infty$-row, this formula defines a topical
operator, also denoted by~$A$. Such operators are called max-plus
operators and operators composition corresponds to the product of
matrices in the max-plus semiring.

Let $(A_n)_{n\in\mathbb{N}}$ be a sequence of random
topical operators on $\mathbb{R}^d$. A stochastic recursive sequence
(SRS) driven by stochastic recursive sequence is a sequence
$(X_n)_{n\in\mathbb{N}}$ satisfying equation $X_{n+1}=A_nX_n$.
To study such sequences, we define $(x(n,x_0))_{n\in\mathbb{N}}$ by
%
%e1 ###
\begin{eqnarray}\label{defx}
x(0,x_0)&=&x_0,
\nonumber\\[-8pt]
\\[-8pt]
x(n+1,x_0)&=&A_nx(n,x_0).
\nonumber
\end{eqnarray}

This class of system can model a wide range of situations. A review
of applications can be found in the last section of~\cite{BM96}. When\vadjust{\goodbreak}
the $x(n,\cdot)$'s are daters, the isotonicity assumption expresses the
causality principle, whereas the additive homogeneity expresses the
possibility to change the origin of time. (See Gunawardena and
Keane~\cite{GunawardenaKeane}, where topical functions were
introduced.) The max-plus case has, for instance, been applied to
model queuing networks (Mairesse~\cite{Mairesse},
Heidergott~\cite{CaractMpQueuNet}), train networks (Heidergott
and De~Vries~\cite{HeidergottDeVriesPubTransNet} and
Braker~\cite{braker}) or job-shop (Cohen, Dubois,
Quadrat and Viot~\cite{cohen85a}). It also
computes the daters of some task resources models (Gaubert and
Mairesse~\cite{gaumair95}) and timed Petri nets including events
graphs (Baccelli \cite{Baccelli}) and 1-bounded Petri nets
(Gaubert and Mairesse~\cite{GaubertMairesseIEEE}). The role of
the max operation is to synchronize different events. For developments
on the max-plus modeling power, see Baccelli, Cohen,
Olsder and Quadrat~\cite{BCOQ} or
Heidergott, Olsder and van~der~Woude~\cite{MpAtWork}.\looseness=1

To clarify things, let us introduce a simple example.
\begin{example}\label{exprod}%ex1.1
Our process assembles two parts. The $n$th time it is done,
it takes time~$a_3(n)$. The parts are prepared separately, which
respectively takes times $a_1(n)$ and $a_2(n)$. Then, they are sent from
the preparation places to the assemblage place, which takes times
$t_1(n)$ and $t_2(n)$ respectively.
Once the assembly place has finished an operation, it asks for new
parts. At that time, if a preparation place has a ready part, it sends
it and starts preparing another one. Otherwise, it finishes the one it
is processing, sends it, and immediately starts preparing another one.
This is summed up in Figure~\ref{BuildSys}.
We denote by $(X_n)_1$ and $(X_n)_2$ the starting date of preparation
of the $n$th part of each type and by $(X_n)_3$ the starting
date of the $(n-1)$th assembly.

%-------------------------------------
%f1 ###
\begin{figure}

\includegraphics{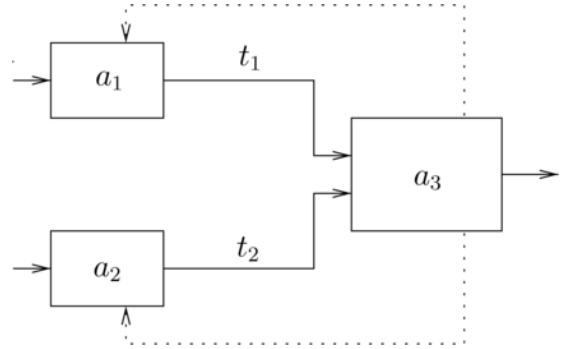}

\caption{A simple production system.\label{BuildSys}}

\end{figure}
%%-------------------------------------

Sequence $(X_n)_{n\in\mathbb{N}}$ is ruled by equations
\begin{eqnarray*}
(X_{n+1})_1 &=& \max\bigl((X_n)_1+a_1(n),(X_n)_3\bigr),
\\ %\qquad
(X_{n})_2 &=& \bigl((X_n)_2+a_1(n),(X_n)_3\bigr)
\end{eqnarray*}
and
\begin{eqnarray*}
(X_{n+1})_3&=& a_3(n)+ \max\bigl((X_n)_1+t_1(n), (X_n)_2+t_2(n)\bigr)
\\
&=& \max \bigl((X_n)_1+t_1(n)+a_3(n-1)+a_1(n),
\\
&&\hspace*{24pt} (X_n)_2+t_2(n)+a_3(n-1)+a_2(n),
\\
&&\hspace*{31pt} (X_n)_3+t_1(n)\vee t_2(n)+a_3(n-1) \bigr),
\end{eqnarray*}
in which we recognize equation~(\ref{defx}), with $A_n$ defined by the
action in the max-plus algebra of
\[
A(n)= \pmatrix{%
a_1(n)&-\infty&a_3(n) \cr
-\infty&a_2(n)&a_3(n) \cr
a_1(n)+t_1(n)&a_2(n)+t_2(n) &t_1(n)\vee t_2(n)+a_3(n-1)}.
\]

We assume that the sequence $(A(n))_{n\in\mathbb{N}}$ is
stationary and ergodic.
\end{example}

We will focus on the asymptotic behavior of $(x(n,\cdot))_{n\in\mathbb{N}}$.
It follows from Theorem~\ref{thVincent}, due to
Vincent, that $(\frac{1}{n}\bigvee_ix_i(n,X^0))_{n\in \mathbb{N}}$
converges to a limit~$\gamma$.

In many cases, if the modeled system is closed, then every sequence of
coordinate $(x_i(n,X^0))_{n\in\mathbb{N}}$ also tends to
$\gamma$, by Theorem~\ref{LGN}. The so-called cycle time $\gamma$
is the inverse of the network's throughput or the inverse of the
production system's output, as in Example~\ref{exprod}. Therefore,
there have been many attempts to estimate it
(Cohen~\cite{Cohen}, Gaujal and Jean-Marie~\cite{ComputIssuesSRS},
Resig et al.~\cite{RVH}). Even when the $A_n$'s are i.i.d. and take only
finitely many values, approximating $\gamma$ is NP-hard
(Blondel, Gaubert and Tsitsiklis~\cite{LyapExpNP}). Hong and his coauthors have obtained
\cite{BacHong1,BacHong2,GaubertHong} analyticity of $\gamma$ as a
function of the law of~$A_1$. They did so under the so-called memory
loss property (MLP) introduced by Mairesse to ensure some stability
of $(x(n,\cdot))_{n\in\mathbb{N}}$ (see~\cite{Mairesse}).

We prove another type of stability under the same assumptions.
If $(A_n)_{n\in\mathbb{N}}$ has the MLP, then $(x(n,\cdot))_{n\in\mathbb{N}}$
actually satisfies a central limit theorem (CLT) under
the same mixing and integrability hypotheses as the real variables in
the CLT for sum of stationary variables by
Billingsley~\cite{ConvProbMeas} and Ibragimov~\cite{Ibrag}.

As far as we know, two CLTs have already been proved for this type of
sequences: one in~\cite{RVH} and the other in~\cite{limtop}. The
most obvious improvement is that both assumed that the $A_n$'s are
i.i.d. Moreover, in~\cite{RVH}, the hypothesis is difficult to check
and the $A_n$'s are max-plus operators defined by almost surely
bounded matrices. In~\cite{limtop} the main hypothesis is also the
MLP, but integrability hypotheses were stronger, except for a subclass
of topical operators.

The remainder of this article is divided into two sections. In
Section~\ref{secPres} we define the memory loss property, present some
law of large number type results and state our central limit theorems.
In Section~\ref{secProofs} we prove the theorems. First, we state the
CLT for subadditive processes by Ishitani~\cite{Ishitani}, then we
check that $(\bigvee_ix_i(n,0))_{n\in\mathbb{N}}$
satisfies each of its hypotheses. To this aim, we use Mairesse's
construction of the stationary version of the SRS, as well as different
results from ergodic theory, depending on the hypothesis. We eventually
deduce the results on $(x(n,\cdot))_{n\in\mathbb{N}}$ from
those on $(\bigvee_ix_i(n,0))_{n\in\mathbb{N}}$.

%s2 ###
\section{Presentation}\label{secPres}
%s2.1 ###
\subsection{Memory loss property}

Dealing with homogeneous operators, it is natural to introduce the
quotient space of $\mathbb{R}^d$ by the\vadjust{\goodbreak} equivalence relation $\sim$
defined by $x\sim y$ if $x-y$ is proportional to
$\mathbf{1}=(1,\ldots,1)'$. This space will be called projective space and
denoted by $\mathbb{PR}_{\max}^d$. Moreover, $\overline{x}$ will be
the equivalence class of~$x$.

The function $\overline{x}\mapsto(x_i-x_j)_{i<j}$ embeds
$\mathbb{PR}_{\max}^d$ onto a subspace of $\mathbb{R}^{{(d(d-1))}/{2}}$
with dimension $d-1$. The infinity norm of $\mathbb{R}^{{(d(d-1))}/{2}}$
therefore induces a distance on $\mathbb{PR}_{\max}^d$
which will be denoted by $\delta$. A direct computation shows that
$\delta(\overline{x},\overline{y})=\bigvee_i(x_i-y_i)+\bigvee_i(y_i-x_i)$.
By a slight abuse, we will also write $\delta(x,y)$ for
$\delta(\overline{x},\overline{y})$. The projective norm of $x$ will
be $|x|_\mathcal{P}=\delta(x,0)=\bigvee_ix_i-\bigwedge_ix_i$.

Let us recall two well-known facts about topical operators. First, a
topical operator is nonexpanding with respect to the infinity norm
(Crandall and Tartar~\cite{CT}). Second, the operator it
defines from $\mathbb{PR}_{\max}^d$ to itself is nonexpanding for
$\delta$ (Mairesse~\cite{Mairesse}).

The key property for our proofs is below:
\begin{defn}[(MLP)]%def2.1
\begin{enumerate}
\item A topical operator $A$ is said to have rank~1 if it defines a
constant operator on $\mathbb{PR}_{\max}^d\dvtx\overline{Ax}$ does
not depend on $x\in\mathbb{R}^d$.
\item A sequence $(A_n)_{n\in\mathbb{N}}$ of
$\mathit{Top}_d$-valued random variables is said to have the memory loss
property (MLP) if there exists an $N$ such that $A_N\cdots A_1$ has
rank~1 with positive probability.
\end{enumerate}
\end{defn}

This notion has been introduced by Mairesse~\cite{Mairesse}, with
the $A_n$'s as max-plus operators. In this case, the denomination
rank~1 is natural.

We have proved in~\cite{GM} that this property is generic for i.i.d.
max-plus operators: it is fulfilled when the support of the law of
$A_1$ is not included in the union of finitely many affine hyperplanes,
and in the discrete case the atoms of the probability measure are
linearly related.

This result applied to Example~\ref{exprod} states that the sequence
$(A_n)_{n\in\mathbb{N}}$ has the MLP provided that the
support of $(a_1(2),a_2(2),t_1(2),t_2(2),a_3(1))$ is not included in a
union of finitely many affine hyperplanes of $\mathbb{R}^5$. This is
not completely straightforward because the matrix $A(1)$ is defined by
only 5 variables, but the detailed result (see Remark~5.1 in \cite{GM})
shows that the linear forms on~$\mathbb{R}^{3\times3}$
that define the hyperplanes are not canceled by $A(1)$, because of the
$-\infty$ entries.

In~\cite{limtop}, we have proved that if the $A_n$'s are i.i.d. and
the sequence has the MLP, then $(x(n,X^0))_{n\in\mathbb{N}}$
satisfies the same limit theorem as a sum of i.i.d. real
variables. Here we prove that it still satisfies the CLT if the
$A_n$'s are mixing quick enough. Quick enough means that the $A_n$'s
satisfy the same integrability and mixing hypothesis as the real
variables in the CLT for the sum of stationary variables by
Billingsley~\cite{ConvProbMeas} and Ibragimov~\cite{Ibrag}.
Moreover, this proves the CLT under weaker integrability condition than
in~\cite{limtop}.

%s2.2 ###
\subsection{Law of large numbers}

There have been many papers about the law of large numbers for products
of random max-plus matrices since its introduction by
Cohen~\cite{Cohen}. We can, for instance, cite Baccelli \cite{Baccelli}, the
most recent paper by Bousch and Mairesse~\cite{BouschMairesseEng}
and Merlet~\cite{theseGM} (in French). The latter article gives
results for a larger class of topical operators, called uniformly topical.

Vincent \cite{vincent} proved a law of large number for topical operators that
will do in our case. He noticed that
$(\bigvee_ix_i(n,0))_{n\in\mathbb{N}}$ [resp.
$(\bigwedge_ix_i(n,0))_{n\in \mathbb{N}}$] is subadditive
(resp. superadditive), which leads to the following:
\begin{thm}[(Vincent~\cite{vincent})]\label{thVincent} %Vincent~thm2.1
Let $(A_n)_{n\in\mathbb{N}}$ be a stationary ergodic
sequence of topical operators and $X^0$ an $\mathbb{R}^d$-valued
random variable. If $A_1.0$ and $X^0$ are integrable, then there exist
$\overline{\gamma}$ and $\underline{\gamma}$ in $\mathbb{R}$ such that
\begin{eqnarray*}
\lim_n \frac{\bigvee_ix_i(n,X^0)}{n}
&=&\overline{\gamma} \qquad \mbox{a.s. and in }\mathbb{L}^1,
\\
\lim_n \frac{\bigwedge_ix_i(n,X^0)}{n}&=&\underline{\gamma}
\qquad \mbox{a.s. and in }\mathbb{L}^1.
\end{eqnarray*}
\end{thm}

Baccelli and Mairesse give a condition to ensure
$\overline{\gamma}=\underline{\gamma}$, hence, the convergence of
$(\frac{x(n,X^0)}{n})_{n\in\mathbb{N}}$:
\begin{thm}[(Baccelli and Mairesse~\cite{BM96})]\label{LGN}%thm2.2
%%Baccelli-Mairesse~
Let $(A_n)_{n\in\mathbb{N}}$ be a stationary ergodic
sequence of topical operators and $X^0$ an $\mathbb{R}^d$-valued
random variable such that $A_1.0$ and $X^0$ are integrable. If there
exists an $N$, such that $A_N\cdots A_1$ has a bounded projective image
with positive probability, then there exists $\gamma$ in $\mathbb{R}$
such that
\[
\lim_n \frac{x(n,X^0)}{n}=\gamma\mathbf{1}
\qquad \mbox{a.s. and in }\mathbb{L}^1.
\]
\end{thm}

That being the case, $\gamma$ is called the Lyapunov exponent of the
sequence. Since matrices with rank~1 have a bounded projective image,
any ergodic sequence $(A_n)_{n\in\mathbb{N}}$ with the
MLP fulfills the hypotheses of Theorem~\ref{LGN}.

%s2.3 ###
\subsection{Statements of the results}

Let us state the definitions of mixing to be used in the sequel.
\begin{defn}[(\textit{Mixing})]\label{defmix}%def2.2
We denote by $\mathcal{F}_n$ the $\sigma$-algebra generated by the
$A_k$'s for $k\le n$ and by $\mathcal{F}^n$ the one generated by the
$A_k$'s for $k\ge n$. We define $\alpha_n$ and $\phi_n$ by the
following:
\begin{enumerate}
\item$\phi(\mathcal{F},\mathcal{G})=\sup\{\frac{|\mathbb{P}(A\cap B)
-\mathbb{P}(A)\mathbb{P}(B)|}{\mathbb{P}(A)}| A\in\mathcal{F}, B\in\mathcal{G}\}$
and $\phi_n=\sup_k \phi(\mathcal{F}_k,\mathcal{F}^{k+n})$.
\item $\alpha(\mathcal{F},\mathcal{G})=\sup\{|\mathbb{P}(A\cap B)
-\mathbb{P}(A)\mathbb{P}(B)|| A\in\mathcal{F},B\in\mathcal{G}\}$
and $\alpha_n=\sup_k \alpha(\mathcal{F}_k,\break\mathcal{F}^{k+n})$.
%%
%Y\|_2}|X\in\mathbb{L}^2(\mathcal{F}), Y\in\mathbb
%{L}^2(\mathcal{G}), \mathbb{E}(X)=\mathbb{E}(Y)=0\}$
%and $\rho_n=\sup_k\rho(\mathcal{F}_k,\mathcal{G}^{k+n})$}
\end{enumerate}
\end{defn}
\begin{thm}\label{thTCL}%thm2.3
If $(A_n)_{n\in\mathbb{N}}$ has the MLP and satisfies one of the
following hypotheses:
\begin{enumerate}[A.]
\item[A.]$A_10 \in\mathbb{L}^2$ and $\sum_{n=1}^\infty\sqrt{\phi_n}<+\infty$,
\item[B.]$A_10 \in\mathbb{L}^{2+\delta}$ and
$\sum_{n=1}^\infty\alpha_n^{{\delta}/{(2+\delta)}} <+\infty$
for some $\delta>0$,
\item[C.]$A_10 \in\mathbb{L}^{\infty}$ and $\sum_{n=1}^\infty\alpha_n<+\infty$,
\end{enumerate}
then
\[
\frac{1}{\sqrt{n}}\bigl(x(n,X^0)-n\gamma\mathbf{1}\bigr)
\stackrel{\mathcal{L}}{\rightarrow} \mathcal{N}\mathbf{1},
\]
where $\mathcal{N}$ is a random variable with zero-mean Gaussian law
(or Dirac measure in~$0$) whose variance does not depend on $X^0$, and
$\stackrel{\mathcal{L}}{\rightarrow}$ denotes the convergence in law.

Moreover, if $X^0$ is integrable, then the variance $\sigma$ of
$\mathcal{N}$ is given by
\[
\lim_{n\rightarrow+\infty}\frac{1}{\sqrt{n}}\mathbb{E}
\biggl|\bigvee_ix_i(n,X^0)-n\gamma\biggr|
= \biggl(\frac{2\sigma^2}{\pi}\biggr)^{1/2}
\]
and $\sigma=0$ if and only if the sequence $(x(n,X^0)-n\gamma
\mathbf{1})_{n\in\mathbb{N}}$ is tight.
\end{thm}
\begin{rem}[(\textit{I.i.d. case})]%rem2.1
When the $A_n$ are i.i.d., I gave more precise results about $\sigma$
in~\cite{limtop}.
In that case, if $\psi$ is a topical function from $\mathbb{R}^d$ to
$\mathbb{R}$, such that $\sup_x|\psi(A_1x)-\psi(x)|$
has a second moment or if $A_1\,0$ has a $(6+\varepsilon)$th moment and
$X^0$ has a $(3+\varepsilon)$th moment, then:
\begin{itemize}
\item$\sigma^2=\lim\frac{1}{n}\mathbb{E}(\psi(x(n,X^0))-n\gamma)^2 $,
\item$\sigma=0$ iff there is a $\theta\in \mathit{Top}_d$ with rank~1 such
that, for any $A$ in the support $S_A$ of $A_1$ and any $\theta'$ with
rank~1 in the semi-group $T_A$ generated by $S_A$, we have
\[
\theta A\theta'=\theta\theta' +\gamma\mathbf{1}.
\]
\end{itemize}
I also proved that if there is such a $\theta$, then every
$\theta\in T_A$ with rank~1 has this property.

Moreover, when the $A_n$ are defined by matrices in the max-plus
algebra, $\sigma$ is positive provided that the support of $A_1$ is
not included in a union of finitely many hyperplanes of
$\mathbb{R}^{d\times d}$.
\end{rem}

In this paper's framework it is not possible to express $\sigma^2$ as
a limit like in the i.i.d. case, because the stationary random
variables in Ishitani's proof of Theorem~\ref{thIshitani} are not
necessarily $\mathbb{L}^2$ (see Section~\ref{secTight}).

%s3 ###
\section{Proofs}\label{secProofs}

%s3.1 ###
\subsection{Results of Ishitani}

We use the results of Ishitani~\cite{Ishitani} for mixing
subadditive processes, which we state now:\vadjust{\goodbreak}

Let $(\Omega,\mathcal{F},T,\mathbb{P})$ be an ergodic measurable
dynamical system, and
$(\mathcal{F}_a^b)_{a,b\in\mathbb{N}}$ a family of sub
$\sigma$-algebras of $\mathcal{F}$, such that
$\mathcal{F}_{a+1}^{b+1}=T^{-1} \mathcal{F}_a^b$, and for any
$a\le c\le d\le b$, $\mathcal{F}_c^d\subset\mathcal{F}_a^b$. The family
$(x_{st})_{s<t}$ of random variables is adapted if, for any $s,t$,
$x_{st}$ is $\mathcal{F}_s^t$-measurable. It is subadditive (resp.
submultiplicative) for any $s<t<u$, $x_{su}\le x_{st}+x_{tu}$
(resp. $x_{su}\le x_{st}\cdot x_{tu}$).
\begin{thm}[(Ishitani \cite{Ishitani} and
Hall and Heyde \cite{HeydeHall})]\label{thIshitani}%thm3.1
Assume $(x_{st})_{s<t}$ is adapted and subadditive. We set $\mathcal
{F}_n=\mathcal{F}_0^n$ and $\mathcal{F}^n=\mathcal{F}_n^{+\infty}$,
and define $\alpha_n$ and $\phi_n$ like in Definition~\ref{defmix}.
We set $(p,\theta)$ as follows:
\begin{enumerate}[(b)]
\item[(a)] $(p,\theta)=(2,2)$ if
$\sum_{n=1}^\infty\sqrt{\phi_n}<+\infty$.
%satisfying $\rho_n\le\sqrt{\phi_n}$. It is not difficult to adapt the
%proof under this weaker hypothesis.}
%
\item[(b)] $(p,\theta)=(2+\delta,\frac{\delta}{2+\delta})$ if
$\sum_{n=1}^\infty\alpha_n^{{\delta}/{(2+\delta)}} <+\infty$
for some $\delta>0$.
\item[(c)] $(p,\theta)=(+\infty,1)$ if $\sum_{n=1}^\infty\alpha_n<+\infty$.
\end{enumerate}

If the following hypotheses are satisfied:
\begin{enumerate}
\item $\lim_{t}\frac{\mathbb{E}(x_{0t})-t\gamma}{\sqrt{t}}=0$,
where $\gamma=\inf\frac{1}{t}\mathbb{E}(x_{0t})$,
\item$\forall t\in\mathbb{N}, |x_{0t}-x_{1t}|\le\Psi$, where
$\Psi\in\mathbb{L}^p$,
\item$\sum_n\sup_t\|x_{0t}-x_{1t}-\mathbb{E}
(x_{0t}-x_{1t}|\mathcal{F}_0^n)\|_\theta<\infty$,
\end{enumerate}
then
\[
\frac{1}{\sqrt{n}}(x_{0n}-n\gamma)
\stackrel{\mathcal{L}}{\rightarrow} \mathcal{N},
\]
where $\mathcal{N}$ is a zero-mean Gaussian law (or a Dirac measure in $0$).

Moreover, the variance $\sigma$ of $\mathcal{N}$ is given by
\[
\lim_{n\rightarrow+\infty}\frac{1}{\sqrt{n}}\mathbb{E}
|x_{0n}-n\gamma| = \biggl(\frac{2\sigma^2}{\pi}\biggr)^{1/2}.
\]
\end{thm}

In the sequel we take $\Omega=\mathit{Top}_d^\mathbb{Z}$, $T$ the shift and
$\mathbb{P}$ such that the law of $(A_n)_{n\in\mathbb{N}}$ is the
image of $\mathbb{P}$ by the projection on the positive coordinates.
From now on, $A_n$~is the projection on the $n$th
coordinate, and we denote $A_0$ by $A$, so that
$A_n=A\circ T^n$.\looseness=1

For any $s<t$, we set $x_{st}=\bigvee_i (A_{t-1}\cdots A_s0)_i$,
and $\mathcal{F}_s^t=\sigma(A_s,\ldots,A_{t-1})$, so
that $(x_{st})_{s,t\in\mathbb{N}}$ is adapted to
$(\mathcal{F}_a^b)_{a,b\in\mathbb{N}}$. Vincent has noticed in~\cite{vincent}
that $(x_{st})_{s,t\in\mathbb{N}}$ is subadditive. From now
on we check that it satisfies hypotheses 1--3 with \mbox{$(p,\theta)=(2,2)$}
under hypothesis~A,
$(p,\theta)=(2+\delta,\frac{\delta}{2+\delta})$ under hypothesis~B
and $(p,\theta)=(+\infty,1)$ under hypothesis~C.

Since $x\mapsto\bigvee_i (A_{t-1}\cdots A_1 x)_i$ is topical,
$\bigwedge_i(A0)_i\mathbf{1} \le A0 \le\bigvee_i(A0)_i\mathbf{1}$ implies
%
%e2 ###
\begin{equation}\label{eqborx}
\bigwedge_i(A0)_i \le x_{0t}-x_{1t}\le\bigvee_i(A0)_i.
\end{equation}
Therefore, we can take $\Psi=|A0|_\infty$ and hypothesis~2 of
Theorem~\ref{thIshitani} is checked.
In the sequel we check the other two hypotheses.

%s3.2 ###
\subsection{Bound on $\mathbb{E}(x_{0t})-t\gamma$}

It is well known and easy to check that, for any $A\in \mathit{Top}_d$ and
$x\in\mathbb{R}^d$, the quantity $\bigvee_i(Ax)_i-\bigvee_ix_i$
only depends on $A$ and $\overline{x}$. We denote it
by $\xi(A,\overline{x})$. With this notation, we have
%
%e3 ###
\begin{equation}\label{eqcob}
\bigvee_ix_i(n,X^0)-\bigvee_i(X^0)_i
= \sum_{k=0}^{n-1}\xi(A_k,\overline{x}(k,X^0)).
\end{equation}
It follows from the main theorem of~\cite{Mairesse}---which can be
extended without difficulty from max-plus to topical operators---that
there is a choice $Y$ of $X^0$, such that
$\overline{x}(n,Y)=\overline{Y}\circ T^n$. In this case,
we see that $\xi(A_k,\overline{x}(k,\overline{Y}))=\xi(A,\overline{Y})\circ T^k$,
therefore, $\bigvee_ix_i(n,Y)-\bigvee_iY_i$
is the partial sum of the stationary
sequence $(\xi(A,\overline{Y})\circ T^k)_{k\in\mathbb{N}}$.

Let us assume for a while that $Y$ is integrable. Then, so is $\xi
(A,\overline{Y})$, because $A0+\bigwedge_iY_i\mathbf{1} \le AY\le A0
+\bigvee_i Y_i\mathbf{1}$ implies
%
%e4 ###
\begin{equation}\label{MajXi}
|\xi(A,\overline{Y})|\le \Biggl|\bigvee_i(A0)_i\Biggr| + |Y|_\mathcal{P}.
\end{equation}
Therefore, it follows from equation~(\ref{eqcob}) with $X^0=Y$ that
\[
\mathbb{E}\Biggl(\bigvee_ix_i(n,Y)\Biggr)
- \mathbb{E}\Biggl(\bigvee_iY_i\Biggr) = n\mathbb{E}(\xi(A,\overline{Y})).
\]
Since topical functions are nonexpanding, we have
$|\bigvee_ix_i(n,Y)-x_{0n}|\le\|Y\|_\infty$, therefore,
$|\mathbb{E}(x_{0n})-n\mathbb{E}(\xi(A,\overline{Y}))|\le2\mathbb{E}(\|Y\|_\infty)$.
In that case $\gamma=\mathbb{E}(\xi(A,\overline{Y}))$
and hypothesis~1 of Theorem~\ref{thIshitani} follows from the
integrability of $Y$.

The end of the subsection is devoted to the proof of the bounds that
will give this integrability.

First, we recall from Mairesse's proof that there is almost surely an
$n\in\mathbb{N}$ such that $rk(A_{-1}\cdots A_{-n})=1$ and
that, for such an $n$, $\overline{Y}=\overline{A_{-1}\cdots A_{-n}0}$.
In the sequel we denote by $N$ the smallest such $n$.

Since $(\bigvee_ix_i(n,0))_{n\in\mathbb{N}}$ [resp.
$(\bigwedge_ix_i(n,0))_{n\in\mathbb{N}}$] is subadditive (resp.
superadditive), we have, for any $n\in\mathbb{N}$ and $i\in[1,d]$,
\[
\sum_{k=1}^{n}\bigwedge_i(A_{-k}0)_i\le(A_{-1}\cdots A_{-n}0)_i
\le \sum_{k=1}^{n}\bigvee_i(A_{-k}0)_i,
\]
therefore, $| A_{-1}\cdots A_{-n}0|_\mathcal{P}\le\sum_{k=1}^n|A_{-k}0|_\mathcal{P}$
and
\begin{eqnarray*}
|Y|_\mathcal{P}&\le& \sum_{n\in\mathbb{N}}\mathbh{1}_{\{N=n\}}
\sum_{k=1}^n|A_{-k}0|_\mathcal{P}
\\
&=&\sum_{k\in\mathbb{N}^*}\mathbh{1}_{\{N\ge k\}}
|A_{-k}0|_\mathcal{P}
\\
&=&\sum_{k\in\mathbb{N}^*}\mathbh{1}_{\{rk(A_{-1}
\cdots A_{-k+1})\neq1\}}|A_{-k}0|_\mathcal{P}.
\end{eqnarray*}
Finally, we get
%
%e5 ###
\begin{equation}\label{eqmajY}
\|Y\|_1\le\|Y\|_\theta\le\sum_{k}
\bigl\|\mathbh{1}_{\{rk(A_{k-1}\cdots A_{1})\neq1\}}|A0|_\mathcal{P}
\bigr\|_\theta.
\end{equation}

The finiteness of the right part of this inequality with $1$ instead of
$\theta$ would be enough to check hypothesis 1, but the finiteness
of this quantity also ensures hypothesis~3, as will be shown in the
next section. Finally, Sections~\ref{secBndA} to~\ref{secBndC} will
be devoted to the proof of the finiteness under each hypothesis of
Theorem~\ref{thTCL}.

%s3.3 ###
\subsection{Bound on $\|x_{0t}-x_{1t}-\mathbb
{E}(x_{0t}-x_{1t}|\mathcal{F}_0^n)\|_\theta$}\label{secBndGene}

We denote by $\Delta_n^t$ the quantity
$|x_{0t}-x_{1t}-\mathbb{E}(x_{0t}-x_{1t}|A_0,\ldots, A_n)|$.
If $t\le n$, then $\Delta_n^t=0$. From now on, we assume $t\ge n$.

First, it follows from equation~(\ref{eqborx}) and the $\mathcal
{F}_0^n$-measurability of $A0$ that
%
%e6 ###
\begin{equation}\label{eqEncadCondit}
\bigwedge_i(A0)_i \le\mathbb{E}(x_{0t}-x_{1t}|
\mathcal{F}_0^n)\le\bigvee_i(A0)_i
\end{equation}
and $\Delta_n^t\le|(A0)|_\mathcal{P}$.

Second, if $rk(A_{n-1}\cdots A_1)=1$, then
%
%e7 ###
\begin{eqnarray}\label{eqAnulRang1}
x_{0t}-x_{1t}-(x_{0n}-x_{1n}) &=& \xi(A_{t-1}\cdots A_{n},
\overline{A_{n-1}\cdots A_00})
\nonumber\\[-8pt]
\\[-8pt]
&&{} -\xi(A_{t-1}\cdots A_{n},
\overline{A_{n-1}\cdots A_10})=0,
\nonumber
\end{eqnarray}
where $\xi$ is the same function as in equation~(\ref{eqcob}).
Therefore, we have
%
%e8 ###
\begin{equation}\label{eqAux1}
\mathbh{1}_{\{rk(A_{n-1}\cdots A_1)= 1\}}(x_{0t}-x_{1t})
=\mathbh{1}_{\{rk(A_{n-1}\cdots A_1)= 1\}}(x_{0n}-x_{1n})
\end{equation}
and
%
%e9 ###
\begin{eqnarray}\label{eqAux2}
\qquad
\mathbh{1}_{\{rk(A_{n-1}\cdots A_1)= 1\}}\mathbb{E}
(x_{0t}-x_{1t}|\mathcal{F}_0^n)
&=& \mathbb{E}\bigl(\mathbh{1}_{\{rk(A_{n-1}\cdots A_1)= 1\}}
(x_{0t}-x_{1t})|\mathcal{F}_0^n\bigr)
\nonumber\\
&=&\mathbb{E}\bigl(\mathbh{1}_{\{rk(A_{n-1}\cdots A_1)= 1\}
}(x_{0n}-x_{1n})|\mathcal{F}_0^n\bigr)
\\
&=&\mathbh{1}_{\{rk(A_{n-1}\cdots A_1)= 1\}}(x_{0n}-x_{1n}).
\nonumber
\end{eqnarray}
Equations~(\ref{eqAux1}) and~(\ref{eqAux2}) together imply that
$\mathbh{1}_{\{rk(A_{n-1}\cdots A_1)= 1\}}\Delta_n^t=0$, and finally,
we have
%
%e10 ###
\begin{equation}\label{eqmajdelt}
\Delta_n^t=\mathbh{1}_{\{rk(A_{n-1}\cdots A_1)\neq1\}}
\Delta_n^t \le\mathbh{1}_{\{rk(A_{n-1}\cdots A_1)\neq1\}}|(A0)|_\mathcal{P}.
\end{equation}

It follows from equations~(\ref{eqmajY}) and~(\ref{eqmajdelt}) that
$(x_{st})$ satisfies hypotheses~1 and~3
of Theorem~\ref{thIshitani}, provided that
%
%e11 ###
\begin{equation}\label{eqsommefini}
\sum_{n=1}^\infty\bigl\|\mathbh{1}_{\{rk(A_{n-1}\cdots A_1)\neq1\}}
|(A0)|_\mathcal{P} \bigr\|_\theta<\infty.
\end{equation}

The next three subsections will prove that relation~(\ref{eqsommefini})
is satisfied, under each of the hypotheses of
Theorem~\ref{thTCL}.

%s3.4 ###
\subsection{Finiteness under hypothesis \textup{A}}\label{secBndA}

From the definition of $\phi$, we see that, for
$X\in\mathbb{L}^1(\mathcal{F})$ and $Y\in\mathbb{L}^\infty(\mathcal{G})$,
$|\mathbb{E}(XY)-\mathbb{E}(X)\mathbb{E}(Y)|\le\phi(\mathcal
{F},\mathcal{G})\|X\|_1\|Y\|_\infty$. We apply this inequality with
$X=|A0|^2_\mathcal{P}$ and $Y=\mathbh{1}_{\{rk(A_{n-1}\cdots
A_{n/2+1})\neq1\}}$, where $n/2$ is the integer part of the half of
$n$, and we take the square root. We get\looseness=1
\begin{eqnarray*}
&& \bigl\|\mathbh{1}_{\{rk(A_{n-1}\cdots A_{1})\neq1\}}
|A0|_\mathcal{P}\bigr\|_{2}
\\
&&\qquad \le \bigl\|\mathbh{1}_{\{rk(A_{n-1}\cdots A_{n/2+1})\neq1\}}
|A0|_\mathcal{P}\bigr\|_{2}
\\
&&\qquad \le \sqrt{\mathbb{P}(rk(A_{n/2}\cdots A_1)\neq1)}
\bigl\||A0|_\mathcal{P}\bigr\|_{2}
+\sqrt{\phi_{n/2}} \bigl\| |A0|_\mathcal{P}\bigr\|_{2}.
\end{eqnarray*}

The $\sqrt{\phi_{n/2}}$'s are summable by hypothesis A. Let us see
that the\break $\sqrt{\mathbb{P}(rk(A_{n/2}\cdots A_1)\neq1)}$'s are too.
For any integers $n$ and $n_0$, we have the following inequality:
\begin{eqnarray*}
&& \mathbb{P}\bigl(rk(A_{n+2n_0}\cdots A_1)\neq1\bigr)
\\
&&\qquad \le \mathbb{E}
\bigl(\mathbh{1}_{\{rk(A_{n}\cdots A_{1})\neq1\}}
\mathbh{1}_{\{rk(A_{n+2n_0}\cdots A_{n+n_0+1})\neq1\}}\bigr)
\\
&&\qquad \le \bigl(\phi_{n_0}+\mathbb{E}\bigl(\mathbh{1}_{\{
rk(A_{n+2n_0}\cdots A_{n+n_0+1})\neq1\}}\bigr)\bigr)
\mathbb{E}\bigl(\mathbh{1}_{\{rk(A_{n}\cdots A_{1})\neq1\}}\bigr)
\\
&&\qquad \le \bigl(\phi_{n_0}+\mathbb{P}
\bigl(rk(A_{n_0}\cdots A_1)\neq 1\bigr)\bigr)
\mathbb{P}\bigl(rk(A_{n}\cdots A_1)\neq1\bigr).
\end{eqnarray*}
Taking $n_0$ big enough, we have
$(\phi_{n_0}+\mathbb{P}(rk(A_{n_0}\cdots A_1)\neq1))<1$, hence,
$\mathbb{P}(rk(A_{n}\cdots A_1)\neq1)$ decreases exponentially fast
and\break
$\sum_{n\in\mathbb{N}}\sqrt{\mathbb{P}(rk(A_{n/2}\cdots A_1)\neq1)}<\infty$.
This concludes the proof of hypotheses 1~and~3
under hypothesis A.

%s3.5 ###
\subsection{Finiteness under hypothesis \textup{B}}\label{secBndB}

Let us take $X=|A0|^{{(2+\delta)}/{(1+\delta)}}_\mathcal{P}$,
$Y=\mathbh{1}_{\{rk(A_{n}\cdots A_{n/2+1})\neq1\}}$ and $q=\frac
{1}{1+\delta}$ in the mixing inequality (see, e.g.,~\cite{HeydeHall})
which states for any $X\in\mathbb{L}^1(\mathcal{F})$ and $Y\in
\mathbb{L}^\infty(\mathcal{G})$
\[
|\mathbb{E}(XY)-\mathbb{E}(X)\mathbb{E}(Y)|\le6\alpha^{1-{1}/{q}}
(\mathcal{F},\mathcal{G})\|X\|_q\|Y\|_\infty\vadjust{\goodbreak}
\]
and let us elevate it to the power $\frac{1+\delta}{2+\delta}$. We get
\begin{eqnarray*}
&& \bigl\|\mathbh{1}_{\{rk(A_{n-1}\cdots A_{1})\neq1\}}
|A0|_\mathcal{P}\bigr\|_{{(2+\delta)}/{(1+\delta)}}
\\
&&\qquad \le \bigl\|\mathbh{1}_{\{rk(A_{n-1}\cdots
A_{n/2+1})\neq1\}}|A0|_\mathcal{P}\bigr\|_{{(2+\delta)}/{(1+\delta)}}
\\
&&\qquad \le \mathbb{P}\bigl(rk(A_{n/2}\cdots A_1)\neq1\bigr)^{{(1+\delta)}/{(2+\delta)}}
\bigl\| |A0|_\mathcal{P}\bigr\|_{{(2+\delta)}/{(1+\delta)}}
\\
&&\qquad \quad {} + (6\alpha_{n/2})^{{\delta}/{(2+\delta)}}
\bigl\||A0|_\mathcal{P}\bigr\|_{2+\delta}.
\end{eqnarray*}

The $(6\alpha_{n/2})^{{\delta}/{(2+\delta)}}$'s are summable by
hypothesis B. To see that the\break $\mathbb{P}(rk(A_{n/2}\cdots
A_1)\neq1)^{{(1+\delta)}/{(2+\delta)}}$
too, we apply the following lemma from~\cite{posmat}
with $\lambda=\frac{2+\delta}{\delta}$
and $M_{st}=\mathbh{1}_{\{rk(A_t\cdots A_s)\neq1\}}$.
\begin{lem}[(Hennion~\cite{posmat})]%lem3.1
Let $(M_{st})_{s<t}$ be submultiplicative and
adapted with values in $[0,1]$ such that $\lim_n\mathbb{E}(M_{0n})=0$.
If $\sum_n\alpha_n^{{1}/{\lambda}}<\infty$, then there exists
$c\in\mathbb{R}$, such that
\[
\mathbb{E}(M_{0n})\le c\biggl(\frac{\ln^2 n}{n}\biggr)^\lambda.
\]
\end{lem}
This concludes the proof of hypotheses 1 and 3 under hypothesis B.

%s3.6 ###
\subsection{Finiteness under hypothesis \textup{C}}\label{secBndC}

We notice that
\begin{eqnarray*}
\sum_{k}\bigl\|\mathbh{1}_{\{rk(A_{k-1}\cdots A_{1})\neq1\}}
|A0|_\mathcal{P}\bigr\|_1
&\le& \sum_{k}\bigl\|\mathbh{1}_{\{rk(A_{k-1}\cdots A_{1})\neq1\}
}\bigr\|_1\bigl\||A0|_\mathcal{P}\bigr\|_\infty
\\
&=& \bigl\| |A0|_\mathcal{P}\bigr\|_\infty\sum_{k}\mathbb{P}(R\ge k)
\\
&=& \bigl\| |A0|_\mathcal{P}\bigr\|_\infty\mathbb{E}(R),
\end{eqnarray*}
where $R=\min\{n|rk(A_{n-1}\cdots A_{1})=1\}$.
Moreover, if $\mathbb{P}(rk(A_{n_0}\cdots A_{1})\neq1)<1$,
then $R-n_0$ is bounded from above by the hitting time of
$\{rk(A_{n_0}\cdots A_{1})= 1\}$. The integrability of $R$ will follow
from the next theorem due to Chazottes.
\begin{thm}[(Chazottes~\cite{ChazottesMixHitt})]\label{thChazottes}%thm3.2
Let $(\Omega,\mathcal{F},\mathbb{P},T)$ be a
measurable dynamical system, $B\in\mathcal{F}$ a set with positive
probability, and $\mathbh{1}_B$ its indicator function. If the mixing
coefficients $\alpha_n$ of the sequence
$(\mathbh{1}_B\circ T^n)_{n\in\mathbb{N}}$ satisfy
$\sum_n\alpha_n<\infty$, then the hitting time of $B$ is integrable.
\end{thm}

To apply the theorem, we notice that, when
$B=\{rk(A_{n_0}\cdots A_{1})= 1\}$, every $\alpha_n$ defined by
$(\mathbh{1}_B\circ T^n)_{n\in\mathbb{N}}$
is less than the $\alpha_{n-n_0}$ defined by
$(A_n)_{n\in\mathbb{N}}$. This ensures the hypothesis of
Theorem~\ref{thChazottes} and concludes the proof of
hypotheses 1 and 3 under hypothesis C.

%s3.7 ###
\subsection{Conclusion of the proof}\label{secConc}

In the last six subsections we have proved that, under hypothesis~A, B
or C of Theorem~\ref{thTCL}, $(\bigvee_ix_i(n,0))_{n\in
\mathbb{N}}$ satisfies the hypotheses of Theorem~\ref{thIshitani}.
Therefore, we have
\[
\frac{1}{\sqrt{n}}\Biggl(\bigvee_ix_i(n,0)-n\gamma\Biggr)
\stackrel{\mathcal{L}}{\rightarrow} \mathcal{N}.
\]

Since topical functions are nonexpanding
%
%e12 ###
\begin{equation}\label{ineqOaXO}
\Biggl|\frac{1}{\sqrt{n}}\bigvee_ix_i(n,X^0)
- \frac{1}{\sqrt{n}}\bigvee_ix_i(n,0)\Biggr|
\le\frac{1}{\sqrt{n}}\|X^0\|_\infty\stackrel{\mathbb{P}}{\rightarrow} 0,
\end{equation}
therefore, $\frac{1}{\sqrt{n}}(\bigvee_ix_i(n,X^0)-n\gamma)
\stackrel{\mathcal{L}}{\rightarrow} \mathcal{N}$.

$\delta_n=(x(n,X^0)-\bigvee_ix_i(n,X^0)\mathbf{1})$ is a
function of $\overline{x}(n,X^0)$, which is converging in law (and
even in total variation) by the main theorem of~\cite{Mairesse},
therefore, \mbox{$\frac{1}{\sqrt{n}}\delta_n\stackrel{\mathbb{P}}{\rightarrow} 0$}
and $\frac{1}{\sqrt{n}}(x(n,X^0)-n\gamma\mathbf{1})
\stackrel{\mathcal{L}}{\rightarrow} \mathcal{N}\mathbf{1}$,
which concludes the proof of the convergence in law.

Inequality (\ref{ineqOaXO}) also implies that
\[
\frac{1}{\sqrt{n}}\mathbb{E}
\Biggl|\bigvee_ix_i(n,X^0)-\bigvee_ix_i(n,0)\Biggr|
\le\frac{1}{\sqrt{n}}\mathbb{E}\|X^0\|_\infty\rightarrow0,
\]
so that
\[
\lim_{n\rightarrow+\infty}\frac{1}{\sqrt{n}}\mathbb{E}
\Biggl|\bigvee_ix_i(n,X^0)-n\gamma\Biggr|
= \biggl(\frac{2\sigma^2}{\pi}\biggr)^{1/2}
\]
follows from $\lim_{n\rightarrow+\infty}
\frac{1}{\sqrt{n}}\mathbb{E}|x_{0n}-n\gamma|
= (\frac{2\sigma^2}{\pi})^{1/2}$.

%s3.8 ###
\subsection{Tightness}\label{secTight}

Without loss of generality, we assume $\gamma=0$.
(Otherwise, just replace $A_n$ by $A_n-\gamma$.)
One part of the equivalence is obvious.
To prove the other part, we have to go into the proof of
Theorem~\ref{thIshitani}. Ishitani constructs a random variable $Z$ (named
$y_{01}$ in~\cite{Ishitani}) and approximates $x_{0n}$ by the Birkhof
sum $S_n=\sum_{k=0}^{n-1}Z\circ T^k$ ($y_{0n}$ in~\cite{Ishitani}).
Then he shows that $(S_n)_{n\in\mathbb{N}}$ fulfills the hypotheses
of Billingsley--Ibragimov's CLT.

In Billingsley--Ibragimov's CLT, the asymptotic variance is zero if and
only if $Z$ is a coboundary, that is, if there is a random variable $f$
such that $Z= f\circ T-f$. (See, e.g., \cite{HeydeHall}.)

Let us assume we are in this situation and identify $Z$. It is built as
a kind of Cesaro type limit of the sequence
$(x_{0n}-x_{1n})_{n\in\mathbb{N}}$. But in our situation
equation~(\ref{eqAnulRang1}) says that this sequence is ultimately constant and that
$X_n$ is equal to the limit as soon as $rk(A_n\cdots A_1)=1$.

Let us denote by $R$ the smallest such $n$ and by $\psi$ the topical
function that maps $x\in\mathbb{R}^d$ to $\bigvee_ix_i$. The random
variable $R$ is almost surely finite because of ergodicity and MLP.
With notation, we have
\[
Z=x_{0R}-x_{1R}=\psi(A_R\cdots A_0 0)-\psi(A_R\cdots A_1 0)
\qquad \mbox{a.s.}
\]
and, for any integer $n$ such that $rk(A_n\cdots A_1)=1$,
%
%e13 ###
\begin{equation}\label{eqZR}
f\circ T-f=\psi(A_n\cdots A_0 0)-\psi(A_n\cdots A_1 0).
\end{equation}

In the sequel we deduce the tightness from equation~(\ref{eqZR}). As
a first and main step, let us show that
$(\psi(A_{R}\cdots A_{-n}0))_{n\in\mathbb{N}}$ is tight.
For any $k\in\mathbb{N}$, since $rk(A_R\cdots A_1)=1$,
$rk(A_R\cdots A_{-k})=1$ and equation~(\ref{eqZR}) holds for
$n=R\circ T^k+k$. Compounded by $T^{-k}$, it becomes
\[
f\circ T^{-k+1}-f\circ T^{-k}=\psi(A_{R}
\cdots A_{-k} 0 )-\psi(A_R\cdots A_{-k+1} 0).
\]

Summing over $k$, we get
\[
f\circ T-f\circ T^{-n}= \psi(A_{R}\cdots A_{-n}0)
- \psi(A_{R}\cdots A_00),
\]
from which the tightness of
$(\psi(A_{R}\cdots A_{-n}0))_{n\in\mathbb{N}}$ is obvious.

The tightness of $(\overline{A_{-1}\cdots A_{-n}0})_{n\in \mathbb{N}}$
is obvious too, because the sequence converges in law.

From those two tightnesses, we successively deduce the tightness of the
following sequences:
\begin{enumerate}[-]
\item[-] $(\psi(A_{-1}\cdots A_{-n}0))_{n\in\mathbb{N}}$,
because equation~(\ref{MajXi}) implies that
\begin{eqnarray*}
|\psi(A_{R}\cdots A_{-n}0)-\psi(A_{-1}\cdots A_{-n}0)|
&=& |\xi(A_{R}\cdots A_{0},\overline{A_{-1}\cdots A_{-n}0})|
\\
&\le& |\psi(A_{R}\cdots A_{0}0)|
+ |\overline{A_{-1}\cdots A_{-n}0}|_\mathcal{P}.
\end{eqnarray*}
\item[-] $((\psi(A_{-1}\cdots A_{-n}0,\overline{A_{-1}
\cdots A_{-n}0})))_{n\in\mathbb{N}}$, again because
$(\overline{A_{-1}\cdots A_{-n}0})_{n\in \mathbb{N}}$ is tight.
\item[-] $(A_{-1}\cdots A_{-n}0)_{n\in\mathbb{N}}$, because
$x\mapsto(\psi(x),\overline{x})$ is a bi-Lipschitz homeomorphism from
$\mathbb{R}^d$ to $\mathbb{R}\times\mathbb{PR}_{\max}^d$
(see~\cite{limtop}).
\item[-] $(x(n,0))_{n\in\mathbb{N}}$, because, for any
$n\in\mathbb{N}$, the random variables $(A_{-1}\cdots A_{-n}0)$
and $x(n,0)$ have the same law.
\item[-] Eventually $(x(n,X^0))_{n\in\mathbb{N}}$, because
the $A_n$ are nonexpanding and, therefore, we have
$\|x(n,X^0)-x(n,0)\|_\infty\le\|X^0\|_\infty$.
\end{enumerate}

%the MLP, the sequence $\bigvee_ix_i(n,0)$ is asymptotically constant,
%because $\bigvee_ix_i(n+1,0)-\bigvee_ix_i(n,0)=\xi(\overline{x}(n,0))$
%Since we have
%%
%A_{-n}0)&=&\xi(A_{R}\cdots A_{0},\overline{A_{-1}\cdots
%A_{-n}0})\\
%&\le& \|A_{R}\cdots A_{0}0\|_\infty+|\overline
%{A_{-1}\cdots A_{-n}0}|_\mathcal{P}
%%
%and $(\overline{A_{-1}\cdots A_{-n}0})_{n\in\mathbb
%{N}}$ is tight, because it converges in law, $(\psi
%(A_{-1}\cdots A_{-n}0))_{n\in\mathbb{N}}$ and
%eventually $((\psi(A_{-1}\cdots A_{-n}0,\overline
%{A_{-1}\cdots A_{-n}0})))_{n\in\mathbb{N}}$ are
%also tight.
%
%Since $x\mapsto(\psi(x),\overline{x})$ is a bilipschitz
%homeomorphism from $\mathbb{R}^d$ to $\mathbb{R}\times\mathbb
%{PR}_{\max}^d$ (see~\cite{limtop}), so is $(A_{-1}\cdots
%A_{-n}0)_{n\in\mathbb{N}}$. Because of the stationarity, the
%variable $(A_{-1}\cdots A_{-n}0)_{n\in\mathbb{N}}$ has
%the same law as $x(n,0)$ for any $n\in N$, therefore $
%(x(n,0))_{n\in\mathbb{N}}$ is tight and $(x(n,X^0)
%)_{n\in\mathbb{N}}$ too, because $\|x(n,X^0)-x(n,0)\|
%_\infty\le\|X^0\|_\infty$.
%}

\section*{Acknowledgments}

This work was done while I was a JSPS post-doctoral fellow at Keio university.
I gratefully thanks H.~Ishitani for introducing me to the article~\cite{Ishitani}
and for interesting discussions, and M.~Keane for his useful rereading.

\printaddresses

\end{document}